\documentstyle[12pt]{article}
\newcommand{\be}{\begin{equation}}
\newcommand{\ee}{\end{equation}}
\newcommand{\bea}{\begin{eqnarray}}
\newcommand{\eea}{\end{eqnarray}}
\newcommand{\bean}{\begin{eqnarray*}}
\newcommand{\eean}{\end{eqnarray*}}

\newcommand{\N}{I\!\!N}
\newcommand{\R}{I\!\!R}
\newcommand{\Z}{Z\!\!\!Z}
\newcommand{\C}{\,I\!\!\!\!C}
\newcommand{\Q}{I\!\!\!\!Q}

\newcommand{\Spec}{\mbox{Spec}\;}

\newcommand{\veee}{\scriptscriptstyle\vee}
\newcommand{\surj}{\longrightarrow\hspace{-1.5em}\longrightarrow}

\newcommand{\ku}{\underline}

\newcommand{\kf}{\footnotesize}
\newcommand{\kff}{\scriptsize}

\newcommand{\ai}{i}

\newcommand{\bj}{j}
\newcommand{\br}{r}
\newcommand{\ci}{i}
\newcommand{\cj}{j}
\newcommand{\kcr}{s}
\newcommand{\dalpha}{d}
\newcommand{\qalpha}{q}

\newcounter{Abschnitt}[section]
\newcommand{\neu}[1]{\protect\refstepcounter{Abschnitt}\protect
   \label{#1}\vspace{1ex}
   {\bf (\protect\arabic{section}.\protect\arabic{Abschnitt})}
                     $\qquad$}
\newcommand{\zitat}[2]{(\protect\ref{#1}.\protect\ref{#1-#2})}

\begin{document}
\title{P-Resolutions of Cyclic Quotients from the Toric Viewpoint}

\author{Klaus Altmann\\
     \small Institut f\"ur reine Mathematik, Humboldt-Universit\"at zu Berlin
     \vspace{-0.7ex}\\ \small Ziegelstra\ss e~13A,
     D-10099 Berlin, Germany. \vspace{-0.7ex}\\ \small E-mail:
     altmann@mathematik.hu-berlin.de}
\date{}
\maketitle

%\tableofcontents
%\par
%\vspace{2ex}

%%%%%%%%%%%%%
%
% 1 Introduction
%
%%%%%%%%%%%%%%
\section{Introduction}\label{I}

%%%%%%%%%%
% (I.1)
%%%%%%%%%

\neu{I-1}
The break through in deformation theory of (two-dimensional)
quotient singularities $Y$ was Koll\'{a}r/Shepherd-Barron's discovery
of the one-to-one correspondence between so-called P-resolutions, on
the one hand, and components of the versal base space, on the other hand
(cf.\ \cite{KS}, Theorem (3.9)).
It generalizes the fact that all deformations admitting
a simultaneous (RDP-) resolution form one single component, the Artin
component.\\
\par

According to defintion (3.8) in \cite{KS}, P-resolutions are partial
resolutions $\pi:\tilde{Y}\to Y$ such that
\begin{itemize}
\item
the canonical divisor $K_{\tilde{Y}|Y}$ is ample relative to $\pi$ (a
minimality condition) and
\item
$\tilde{Y}$ contains only mild singularities of a certain type (so-called
T-singularities).
\end{itemize}

Despite their definition as those quotient singularities admitting a
$\Q$-Gorenstein one-parameter smoothing (\cite{KS}, (3.7)),
there are at least three further descriptions of the class of T-singularities:
An explicit list of their defining group actions on $\C^2$ (\cite{KS}, (3.10)),
an inductive procedure to construct their resolution graphs
(\cite{KS}, (3.11)), and a characterization using toric language
(\cite{Homog}, (7.3)).\\
The latter one starts with the observation that
affine, two-dimensional toric varieties (given by some rational, polyhedral
cone $\sigma\subseteq\R^2$) provide exactly the two-dimensional cyclic
quotient singularties. Then, T-singularities come from cones over rational
intervals of integer length placed in height one (i.e.\ contained in the
affine line $(\bullet,1)\subseteq\R^2$).
\vspace{-2ex}
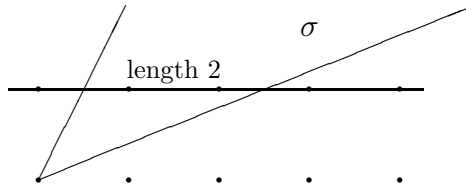
\begin{center}
\unitlength=0.4mm
\linethickness{0.4pt}
\begin{picture}(155.00,68.00)
\put(0.00,40.00){\line(1,0){138.00}}
\put(10.00,40.00){\circle*{2.00}}
\put(40.00,40.00){\circle*{2.00}}
\put(70.00,40.00){\circle*{2.00}}
\put(100.00,40.00){\circle*{0.00}}
\put(100.00,40.00){\circle*{2.00}}
\put(10.00,10.00){\line(5,2){145.00}}
\put(10.00,10.00){\line(1,2){29.00}}
\put(130.00,40.00){\circle*{2.00}}
\put(130.00,10.00){\circle*{2.00}}
\put(100.00,10.00){\circle*{2.00}}
\put(70.00,10.00){\circle*{2.00}}
\put(40.00,10.00){\circle*{2.00}}
\put(10.00,10.00){\circle*{2.00}}
\put(55.00,44.00){\makebox(0,0)[cb]{\kf length 2}}
\put(100,60){\makebox(0,0)[cc]{$\sigma$}}
\end{picture}
\end{center}
If the affine interval is of length $\mu+1$, then the corresponding
T-singularity
will have Milnor number $\mu$ (on the $\Q$-Gorenstein one-parameter
smoothing).\\
\par

%%%%%%%%%%
% (I.2)
%%%%%%%%%

\neu{I-2}
In \cite{Ch-CQS} and \cite{St-CQS} Christophersen and Stevens gave a
combinatorial
description of all P-resolutions for two-dimensional, cyclic quotient
singularities.
Using an inductive construction method (going through different cyclic
quotients
with step-by-step increasing multiplicity) they have shown that there is a
one-to-one
correspondence between P-resolutions, on the one hand, and certain integer
tuples
$(k_2,\dots,k_{e-1})$ yielding zero if expanded as a (negative) continued
fraction
(cf.\ \zitat{P}{2}), on the other hand.\\
\par

The aim of the present paper is to provide an elementary, direct method for
constructing the
P-resolutions of a cyclic quotient singularity (i.e.\ a two-dimensional toric
variety)
$Y_{\sigma}$. Given a chain $(k_2,\dots,k_{e-1})$ representing zero, we will
give a straight
description of the corresponding polyhedral subdivision of $\sigma$. (In
particular,
the bijection between those
0-chains and P-resolutions will be proved again by a different
method.)\\
\par

%%%%%%%%%%%
%
% 2 Cyclic Quotient Singularities
%
%%%%%%%%%%%%%%
\section{Cyclic Quotient Singularities}\label{CQS}

In the following we want to remind the reader of basic notions concerning
continued fractions and cyclic quotients. It should be considered a good
chance to fix notations. References are \cite{Oda} (\S 1.6) or
the first sections in \cite{Ch-CQS} and \cite{St-CQS}, respectively.\\
\par

%%%%%%%%%%
% (CQS.1)
%%%%%%%%%

\neu{CQS-1}
{\bf Definition:}
To integers $c_1,\dots,c_r\in\Z$ we will assign the continued fraction
$[c_1,\dots,c_r]\in\Q$, if the following inductive procedure makes sense
(i.e.\ if no division by 0 occurs):
\begin{itemize}
\item
$[c_r]:=c_r$
\item
$[c_i,\dots,c_r]:= c_i- 1/[c_{i+1},\dots,c_r]$.
\end{itemize}
If $c_i\geq 2$ for $i=1,\dots,r$, then $[c_1,\dots,c_r]$ is always defined
and yields a rational number greater than 1. Moreover, all these numbers may be
represented by those continued fractions in a unique way.\\
\par

%%%%%%%%%%
% (CQS.2)
%%%%%%%%%

\neu{CQS-2}
Let $n\geq 2$ be an integer and $q\in (\Z/n\Z)^\ast$ be represented by an
integer
between $0$ and $n$. These data provide a group action of $\Z/n\Z$ on $\C^2$
via the matrix
$\left(\begin{array}{cc}\xi&0\\0&\xi^q\end{array}\right)$
(with $\xi$ a primitive $n$-th
root of unity). The quotient is denoted by $Y(n,q)$.\\
\par

In toric language, $Y(n,q)$ equals the variety $Y_\sigma$ assigned to the
polyhedral cone
$\sigma:=\langle (1,0);(-q,n)\rangle\subseteq\R^2$. ($Y_\sigma$ is defined as
$\Spec \C[\sigma^{\veee}\cap\Z^2]$ with
\[
\sigma^{\veee}:=\{r\in(\R^2)^\ast\,| \;r\geq 0 \mbox{ on } \sigma\}=
\langle [0,1];[n,q]\rangle \subseteq (\R^2)^\ast \cong \R^2\,.)
\]
{\bf Notation:} Just to distinguish between $\R^2$ and its dual
$(\R^2)^\ast\cong\R^2$,
we will denote these vector spaces by $N_{\R}$ and $M_{\R}$, respectively.
(Hence,
$\sigma\subseteq N_{\R}$ and $\sigma^{\veee}\subseteq M_{\R}$.)
Elements of $N_{\R}\cong\R^2$ are written in paranthesis; elements of
$M_{\R}\cong\R^2$ are written in brackets. The natural pairing
between $N_{\R}$ and $M_{\R}$ is denoted by $\langle\,,\,\rangle$ which should
not be
mixed up with the symbol indicating the generators of a cone. Finally, all
these
remarks apply for the lattices $N\cong\Z^2$ and $M\cong\Z^2$, too.\\
\par

%%%%%%%%%%
% (CQS.3)
%%%%%%%%%

\neu{CQS-3}
Let $n,q$ as before. We may write $n/(n-q)$ and $n/q$ (both are greater than 1)
as continued fractions
\[
n/(n-q) = [a_2,\dots,a_{e-1}] \; \mbox{ and } \; n/q = [b_1,\dots,b_\br]
\quad (a_\ai, b_\bj\geq 2).
\]
The $a_\ai$'s and the $b_\bj$'s are mutually related by Riemenschneider's point
diagram (cf.\ \cite{Riem}).\\
\par

Take the convex hull of $(\sigma^{\veee}\cap M)\setminus\{0\}$ and denote by
$w^1,w^2,\dots, w^e$ the lattice points on its compact edges. If ordered the
right way, we obtain $w^1=[0,1]$ and $w^e=[n,q]$ for the first and the last
point,
respectively.
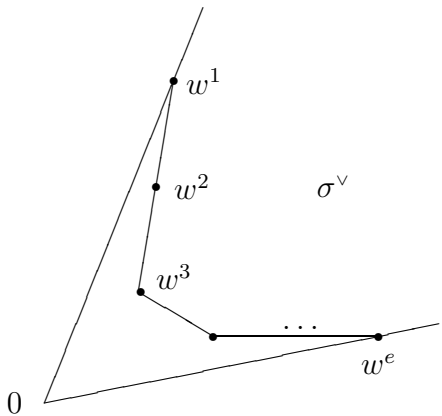
\begin{center}
\unitlength=0.4mm
\linethickness{0.4pt}
\begin{picture}(141.00,142.00)
\put(0,10){\makebox(0,0)[cc]{$0$}}
\put(10.00,10.00){\line(2,5){52.67}}
\put(10.00,10.00){\line(5,1){131.00}}
\put(53.00,118.00){\line(-1,-6){11.67}}
\put(41.33,47.00){\line(5,-3){24.67}}
\put(66.00,32.33){\line(1,0){54.00}}
\put(53.00,117.00){\circle*{3.00}}
\put(42.00,47.00){\circle*{3.00}}
\put(47.00,82.00){\circle*{3.00}}
\put(66.00,32.00){\circle*{3.00}}
\put(121.00,32.00){\circle*{3.00}}
\put(63.00,117.00){\makebox(0,0)[cc]{$w^1$}}
\put(59.00,82.00){\makebox(0,0)[cc]{$w^2$}}
\put(95.00,35.00){\makebox(0,0)[cc]{$\dots$}}
\put(53.00,52.00){\makebox(0,0)[cc]{$w^3$}}
\put(121.00,23.00){\makebox(0,0)[cc]{$w^e$}}
\put(106.00,82.00){\makebox(0,0)[cc]{$\sigma^{\veee}$}}
\end{picture}
\end{center}
Then, $E:=\{w^1,\dots,w^e\}$ is the minimal generating set (the so-called
Hilbert basis) of the semigroup $\sigma^{\veee}\cap M$. These point are
related to our first continued fraction by
\[
w^{\ai-1}+w^{\ai+1}=a_\ai\,w^\ai \quad (\ai=2,\dots,e-1).
\]
{\bf Remark:}
The surjection $\N^E\surj \sigma^{\veee}\cap M$ provides a minimal
embedding of $Y_{\sigma}$. In particular, $e$ equals its embedding
dimension.\\
\par

In a similar manner we can define $v^0,\dots,v^{\br+1}\in\sigma\cap N$
in the original cone; now we
have $v^0=(1,0)$, $v^{\br+1}=(-q,n)$, and the relation to the continued
fractions is $v^{\bj-1}+v^{\bj+1}=b_\bj\,v^\bj$ (for $\bj=1,\dots,\br$).
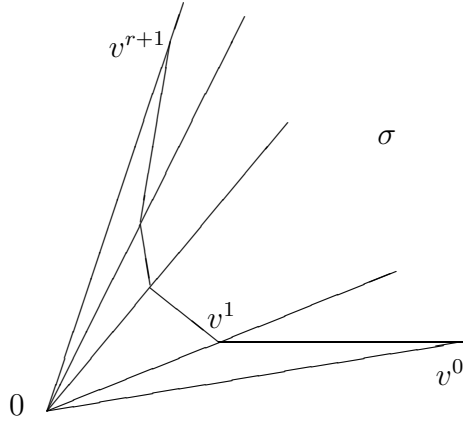
\begin{center}
\unitlength=0.4mm
\linethickness{0.4pt}
\begin{picture}(148.00,146.00)
\put(10.00,10.00){\line(1,3){45.33}}
\put(10.00,10.00){\line(6,1){138.00}}
\put(10.00,10.00){\line(1,2){65.67}}
\put(10.00,10.00){\line(5,6){80.00}}
\put(10.00,10.00){\line(5,2){116.00}}
\put(51.00,133.00){\line(-1,-6){10.00}}
\put(41.00,72.00){\line(1,-6){3.33}}
\put(44.33,51.00){\line(5,-4){22.67}}
\put(67.00,33.00){\line(1,0){81.00}}
\put(144.00,23.00){\makebox(0,0)[cc]{$v^0$}}
\put(68.00,42.00){\makebox(0,0)[cc]{$v^1$}}
\put(40.00,132.00){\makebox(0,0)[cc]{$v^{\br+1}$}}
\put(0.00,12.00){\makebox(0,0)[cc]{$0$}}
\put(123.00,101.00){\makebox(0,0)[cc]{$\sigma$}}
\end{picture}
\end{center}
Drawing rays through the origin and each point $v^\bj$, respectively,
provides a polyhedral subdivision $\Sigma$ of $\sigma$. The corresponding
toric variety $Y_\Sigma$ is a resolution of our singularity
$Y_\sigma$. The numbers $-b_\bj$ equal the self intersection numbers of the
exceptional divisors; since $b_\bj\geq 2$, the resolution is the {\em minimal}
one.\\
\par

%%%%%%%%%%%%%
%
% 3 The Maximal Resolution
%
%%%%%%%%%%%%%%
\section{The Maximal Resolution}\label{MR}

%%%%%%%%%%
% (MR.1)
%%%%%%%%%

\neu{MR-1}
{\bf Definition:} (\cite{KS}, (3.12)) For a resolution $\pi:\tilde{Y}\to Y$
we may write $K_{\tilde{Y}|Y}:=K_{\tilde{Y}}-\pi^\ast K_Y=
\sum_{\cj}(\alpha_\cj-1)E_\cj$
($E_\cj$ denote the exceptional divisors, $\alpha_\cj\in\Q$).
Then, $\pi$ will be called {\em maximal}, if it is maximal with respect to the
property $0<\alpha_\cj<1$.\\
\par

The maximal resolution is uniquely determined and dominates all the
P-resolutions.
Hence, for our purpose, it is
more important than the minimal one. It can be constructed from the minimal
resolution
by sucsessive blowing
up of points $E_\ci\cap E_\cj$ with $\alpha_\ci+\alpha_\cj\geq 0$ (cf.\ Lemma
(3.13)
and Lemma (3.14) in \cite{KS}).\\
\par

%%%%%%%%%%
% (MR.2)
%%%%%%%%%

\neu{MR-2}
{\bf Proposition:}
{\em
The maximal resolution of $Y_\sigma$ is toric. It can be obtained by drawing
rays through $0$ and all interior lattice points (i.e.\ $\in N$) of the
triangle
$\Delta:= \mbox{\em conv} \,(0,v^0,v^{\br+1})$, respectively.
}\\
\par

{\bf Proof:}
We have to keep track of the rational numbers $\alpha_\cj$. Hence, we will show
how
they can be ``seen'' in an arbitrary toric resolution of $Y_\sigma$. Let
$\Sigma<\sigma$
be a subdivision generated by one-dimensional rays through the points
$u^0,\dots,u^{\kcr+1}\in\sigma\cap N$. (In particular, $u^0=v^0=(1,0)$ and
$u^{\kcr+1}=v^{\br+1}=(-q,n)$; moreover, for the minimal resolution we would
have
$\kcr=\br$ and $u^\bj=v^\bj$ ($\bj=0,\dots,\br+1$).) Denote by
$c_1,\dots,c_\kcr$ the
integers given by the relations
\[
u^{\cj-1}+u^{\cj+1}=c_\cj\,u^\cj \qquad (\cj=1,\dots,\kcr).
\]
(In particular, $c_\bj =b_\bj$ for the minimal resolution again.)\\
\par

As usual, the numbers $-c_\cj$ equal the self intersection numbers of the
exceptional
divisors $E_\cj$ in $Y_\Sigma$: Indeed, $D:=\sum_\ci u^\ci\,E_\ci$ is a
principal divisor
(if you do not like coefficients $u^\ci$ from $N$, evaluate them by arbitrary
elements of $M$); hence,
\[
\begin{array}{rcl}
0 \,=\, E_\cj\cdot D
&=& E_\cj\cdot (u^{\cj-1} E_{\cj-1} + u^\cj E_\cj +
u^{\cj+1} E_{\cj+1})\\
&=& u^{\cj-1}+ (E_\cj)^2\,u^\cj + u^{\cj+1}\\
&=& ( c_\cj + (E_\cj)^2 )\cdot u^\cj \qquad\qquad (\cj=1,\dots,\kcr).
\end{array}
\]
\par

On the other hand, we can use the projection formula to obtain
\[
\begin{array}{rcl}
-2 \,= \, 2\,g(E_\cj)-2
&=& K_{\tilde{Y}|Y}\cdot E_\cj + (E_\cj)^2\\
&=& \sum_\ci \,(\alpha_\ci-1) \,(E_\ci\cdot E_\cj) + (E_\cj)^2\\
&=& (\alpha_{\cj-1}-1) + (\alpha_\cj -1)\,(E_\cj)^2 + (\alpha_{\cj+1}-1) +
(E_\cj)^2\,,
\vspace{-1ex}
\end{array}
\]
% hence $ \alpha_{\cj-1} + \alpha_{\cj+1} = c_\cj\, \alpha_\cj$ with
% $\cj=1,\dots,\kcr$ and $\alpha_0,\,\alpha_{\kcr+1}:=1$.
hence
\[
\alpha_{\cj-1} + \alpha_{\cj+1} = c_\cj\, \alpha_\cj
\qquad (\cj=1,\dots,\kcr\,;\; \alpha_0,\,\alpha_{\kcr+1}:=1).
\vspace{1ex}
\]
\par

Looking at the definition of the $c_\cj$'s (via relations among the lattice
points
$u^\cj$), there has to be some $R\in M_{\R}$ such that
\[
\alpha_\cj = \langle u^\cj,\,R\rangle \qquad\qquad
(\cj=0,\dots,\kcr+1).
\]
The conditions $\langle u^0,\,R\rangle = \alpha_0 =1$ and
$\langle u^{\kcr+1},\,R\rangle = \alpha_{\kcr+1} =1$ determine $R$ uniquely.
Now, we can see that $\alpha_\cj$ measures exactly the quotient between the
length of the
line
segment $\overline{0\,u^\cj}$, on the one hand, and the length of the
$\Delta$-part
of the line through $0$ and $u^\cj$, on the other hand. In particular,
$\alpha_\cj<1$ if and only if $u^\cj$ sits below the line connecting $u^0$ and
$u^{\kcr+1}$.
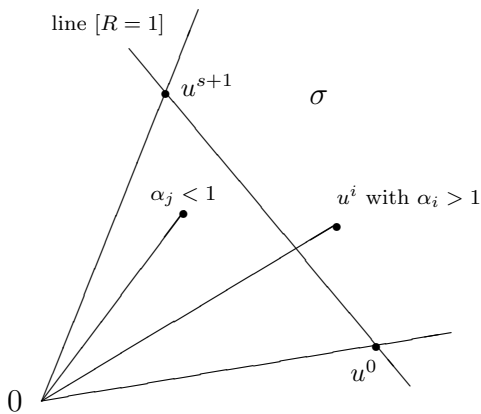
\begin{center}
\unitlength=0.4mm
\linethickness{0.4pt}
\begin{picture}(146.00,140.00)
\put(10.00,10.00){\line(2,5){52.00}}
\put(10.00,10.00){\line(6,1){136.00}}
\put(39.00,127.00){\line(5,-6){93.33}}
\put(10.00,10.00){\line(3,4){47.33}}
\put(10.00,10.00){\line(5,3){98.00}}
\put(51.00,112.00){\circle*{3.00}}
\put(121.00,28.00){\circle*{3.00}}
\put(57.00,72.00){\circle*{3.00}}
\put(108.00,68.00){\circle*{3.00}}
\put(1.00,10.00){\makebox(0,0)[cc]{$0$}}
\put(117.00,20.00){\makebox(0,0)[cc]{\kf $u^0$}}
\put(108.00,78.00){\makebox(0,0)[lc]{\kff $u^\ci$ with $\alpha_\ci>1$}}
\put(57.00,78.00){\makebox(0,0)[cc]{\kff $\alpha_\cj<1$}}
\put(65.00,115.00){\makebox(0,0)[cc]{\kf $u^{\kcr+1}$}}
\put(32.00,135.00){\makebox(0,0)[cc]{\kff line $[R=1]$}}
\put(102.00,111.00){\makebox(0,0)[cc]{$\sigma$}}
\end{picture}
\end{center}
This explains how to construct the maximal resolution: Start with the minimal
one
and continue subdividing each small cone $\langle u^\cj,u^{\cj+1}\rangle$
into $\langle u^\cj, u^\cj+ u^{\cj+1}\rangle \,\cup \,
\langle u^\cj + u^{\cj+1}, u^{\cj+1}\rangle$ as long as it contains interior
lattice
points below the line $[R=1]$, i.e.\ belonging to $\mbox{int}\, \Delta$.
\hfill$\Box$\\
\par

{\bf Corollary:}
{\em Every P-resolution is toric.}\\
\par

{\bf Proof:} P-resolutions are obtained by blowing down curves in the maximal
resolution.
\hfill$\Box$\\
\par

%%%%%%%%%%
% (MR.3)
%%%%%%%%%

\neu{MR-3}
{\bf Example:}
We take the example $Y(19,7)$ from \cite{KS}, (3.15). Since
$\sigma=\langle (1,0),\, (-7,19) \rangle$, the interior of $\Delta$ is given by
the
three inequalities
\[
y>0\,,\; 19x+7y>0\,, \mbox{ and}\; 19x+8y<19\;\,
(\mbox{corresponding to } R=[1,\,8/19])\,.
\]
The only primitive (i.e.\ generating rays) lattice points contained in
$\mbox{int}\,\Delta$ are
% $\,u^1=(0,1)$, $\,u^2=(-1,4)$, $\,u^3=(-2,7)$, $\,u^4=(-1,3)$,
% $\,u^5=(-5,14)$, and $\,u^6=(-4,11)$.
\[
u^1\!=(0,1)\,,\; u^2\!=(-1,4)\,,\; u^3\!=(-2,7)\,,\; u^4\!=(-1,3)\,,\;
u^5\!=(-5,14)\,,\; u^6\!=(-4,11)\,.
\]
They provide the maximal resolution.
The corresponding $\alpha$'s can be obtained
by taking the scalar product with $R=[1,\,8/19]$, i.e.\
they are $8/19$, $13/19$, $18/19$, $5/19$, $17/19$, and $12/19$.\\
The minimal resolution uses only the rays through $\,u^1=(0,1)$,
$\,u^4=(-1,3)$, and
$\,u^6=(-4,11)$, respectively.\\
\par

%%%%%%%%%%%%%
%
% 4 P-Resolutions
%
%%%%%%%%%%%%%%
\section{P-Resolutions}\label{P}

%%%%%%%%%%
% (P.1)
%%%%%%%%%

\neu{P-1}
In this section we will speak about {\em partial} toric resolutions
$\pi:Y_\Sigma\to Y_\sigma$. Nevertheless, we use the same notation as we
did for the maximal resolution: The fan $\Sigma$ subdividing $\sigma$ is
genarated by rays through $u^0,\dots,u^\kcr\in\sigma\cap N$; each ray
$u^\cj$ corresponds to an exceptional divisor $E_\cj\subseteq Y_\Sigma$.
However, since $u^{\cj-1} + u^{\cj+1}$ need not to be a multiple of
$u^\cj$, the numbers $c_\cj$ do not make sense anymore.\\
\par

{\bf Lemma:} (\cite{toricMori}, (4.3))
{\em
For $K:=K_{Y_\Sigma}$ or $K:=K_{Y_\Sigma|Y_\sigma}$ the intersection number
$(E_\cj\cdot K)$ is positive, zero, or negative, if the line segments
$\overline{u^{\cj-1} u^\cj}$ and $\overline{u^\cj u^{\cj+1}}$ form a strict
concave, flat, or strict convex ``roof'' over the two cones, respectively.
}
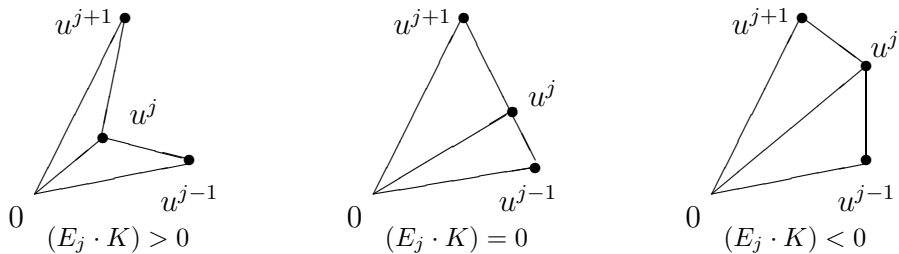
\begin{center}
\unitlength=0.50mm
\linethickness{0.4pt}
\begin{picture}(231.00,68.50)
\put(18.00,66.00){\makebox(0,0)[cc]{$u^{\cj+1}$}}
\put(46.00,19.00){\makebox(0,0)[cc]{$u^{\cj-1}$}}
\put(27.00,8.00){\makebox(0,0)[cc]{\kf $(E_\cj\cdot K)>0$}}
\put(108.00,66.00){\makebox(0,0)[cc]{$u^{\cj+1}$}}
\put(136.00,19.00){\makebox(0,0)[cc]{$u^{\cj-1}$}}
\put(117.00,8.00){\makebox(0,0)[cc]{\kf $(E_\cj\cdot K)=0$}}
\put(198.00,66.00){\makebox(0,0)[cc]{$u^{\cj+1}$}}
\put(226.00,19.00){\makebox(0,0)[cc]{$u^{\cj-1}$}}
\put(180.00,14.00){\makebox(0,0)[cc]{$0$}}
\put(207.00,8.00){\makebox(0,0)[cc]{\kf $(E_\cj\cdot K)<0$}}
\put(231.00,60.00){\makebox(0,0)[cc]{$u^\cj$}}
\put(140.00,47.00){\makebox(0,0)[cc]{$u^\cj$}}
\put(34.00,42.00){\makebox(0,0)[cc]{$u^\cj$}}
\put(0.00,14.00){\makebox(0,0)[cc]{$0$}}
\put(90.00,14.00){\makebox(0,0)[cc]{$0$}}
\put(226.00,29.00){\circle*{3.00}}
\put(226.00,54.00){\circle*{3.00}}
\put(29.00,67.00){\circle*{3.00}}
\put(119.00,67.00){\circle*{3.00}}
\put(138.00,27.00){\circle*{3.00}}
\put(132.00,42.00){\circle*{3.00}}
\put(23.00,35.00){\circle*{3.00}}
\put(209.00,67.00){\circle*{3.00}}
\put(46.00,29.00){\circle*{3.00}}
\put(5.00,20.00){\line(1,2){23.67}}
\put(5.00,20.00){\line(6,5){18.00}}
\put(5.00,20.00){\line(5,1){41.00}}
\put(29.00,67.00){\line(-1,-5){6.33}}
\put(23.00,35.00){\line(4,-1){23.00}}
\put(95.00,20.00){\line(1,2){23.67}}
\put(185.00,20.00){\line(1,2){23.67}}
\put(185.00,20.00){\line(5,1){41.00}}
\put(209.00,67.00){\line(4,-3){17.00}}
\put(226.00,54.33){\line(0,-1){25.33}}
\put(185.00,20.00){\line(6,5){41.00}}
\put(119.00,67.00){\line(1,-2){19.00}}
\put(95.00,20.00){\line(6,1){44.00}}
\put(95.00,20.00){\line(5,3){38.00}}
\end{picture}
\end{center}

{\bf Proof:}
Using $K:=K_{Y_\Sigma}= -\sum_{\ci=0}^{\kcr+1}E_\ci$
(cf.\ \cite{Oda}, (2.1)) we have
\[
(E_\cj\cdot K) = - (E_\cj\cdot E_{\cj-1}) -(E_\cj)^2 -
(E_\cj\cdot E_{\cj+1})\,.
\]
On the other hand, as in the proof of Proposition \zitat{MR}{2}, we know
that
\[
0 = (E_\cj\cdot E_{\cj-1}) \,u^{\cj-1} + (E_\cj)^2 \, u^\cj
+ (E_\cj\cdot E_{\cj+1}) \, u^{\cj+1}\,.
\]
Combining both formulas yields the final result
\[
(E_\cj\cdot K)\, u^\cj =
(E_\cj\cdot E_{\cj-1}) \,(u^{\cj-1}-u^\cj) +
(E_\cj\cdot E_{\cj+1}) \,(u^{\cj+1}-u^\cj)\,.
\vspace{-3ex}
\]
\hspace*{\fill}$\Box$\\
\par

{\bf Remark:} The previous lemma together with Proposition \zitat{MR}{2}
illustrate again the fact that all P-resolutions (and we just need the
fact that the canonical divisor is relatively ample) are dominated by the
maximal resolution.\\
\par

%%%%%%%%%%
% (P.2)
%%%%%%%%%

\neu{P-2}
In \cite{Ch-CQS} Christophersen has defined the set
\[
K_{e-2}:=\{(k_2,\dots,k_{e-1})\in \N^{e-2}\,|\; [k_2,\dots,k_{e-1}]
\mbox{ is well defined and yields } 0\,\}
\]
of chains representing zero. To every such chain there are assigned
non-negative
integers $\qalpha_1,\dots,\qalpha_e$ characterized by the following mutually
equivalent properties:
\begin{itemize}
\item
$\qalpha_1=0$, $\,\qalpha_2=1$, and
$\;\qalpha_{\ai-1} + \qalpha_{\ai+1} = k_\ai\, \qalpha_\ai
\quad (i=2,\dots,e-1)$;
\item
$\qalpha_{e-1}=1$, $\,\qalpha_e=0$, and
$\;\qalpha_{\ai-1} + \qalpha_{\ai+1} = k_\ai\, \qalpha_\ai
\quad (i=2,\dots,e-1)$;
\item
$\qalpha_e=0$ and
$\,[k_\ai,\dots,k_{e-1}]= \qalpha_{\ai-1}/\qalpha_\ai$ with
$\mbox{gcd}(\qalpha_{\ai-1},\qalpha_\ai)=1 \;(i=2,\dots,e-1)$.
% \vspace{1ex}
\end{itemize}
(The two latter properties do not even use the fact that the continued
fraction $[k_2,\dots,k_{e-1}]$ yields zero.)\\
\par

{\bf Remark:}
The elements of $K_{e-2}$ correspond one-to-one to triangulations of a
(regular)
$(e-1)$-gon with vertices $P_2,\dots,P_{e-1},P_\ast$. Then, the numbers $k_\ai$
tell
how many triangles are attached to $P_\ai$. The numbers $\qalpha_\ai$ have an
easy meaning in this language, too.\\
\par

Finally, for a given $Y_\sigma$ with embedding dimension $e$, Christophersen
defines
\[
K(Y_\sigma):=\{(k_2,\dots,k_{e-1})\in K_{e-2}\,|\; k_\ai\leq a_\ai\}\,.
\vspace{1ex}
\]
\par

{\bf Theorem:}
{\em
Each P-resolution of $Y_\sigma$ (i.e.\ the corresponding subdivision $\Sigma$
of $\sigma$)
is given by some $\ku{k}\in K(Y_\sigma)$ in the following way:
\vspace{0.5ex}\\
(1) $\Sigma$ is built from the rays that are orthogonal to
$w^\ai/\qalpha_\ai - w^{\ai-1}/\qalpha_{\ai-1}\in M_{\R}$
(for $\ai=3,\dots,e-1$).
In some sense, if the occuring divisions by zero are interpreted well,
$\Sigma$ may be seen as dual to the Newton boundary generated by
$w^\ai/\qalpha_\ai\in\sigma^{\veee}$ ($\ai=1,\dots,e$).
\vspace{0.5ex}\\
(2) The affine lines $[\langle \bullet,\, w^\ai\rangle = \qalpha_\ai]$ form the
``roofs'' of the $\Sigma$-cones. In particular, the (possibly degenerate)
cones $\tau^\ai\in\Sigma$ correspond to the elements $w^1,\dots,w^e\in E$
% even in their order.
% The generators of $\tau^\ai$ are obtained as the intersection points
% of the actual ``roof'' line with the adjacent ones.
% \vspace{0.5ex}\\
% (3)
The ``roof'' over the cone $\tau^\ai$ has length
$\ell_\ai:= (a_\ai-k_\ai)\,\qalpha_\ai$
(the lattice structure $M\subseteq M_{\R}$ induces a metric on rational lines).
In particular, $\tau^\ai$ is degenerated if
and only if $k_\ai=a_\ai$. The Milnor number of the T-singularity
$Y_{\tau^\ai}$
equals $(a_\ai-k_\ai-1)$.
}\\
\par

%%%%%%%%%%
% (P.3)
%%%%%%%%%

\neu{P-3}
{\bf Proof:}
According to the notation introduced in \zitat{P}{1}, the fan $\Sigma$ consists
of (non-degenerate) cones $\tau^\cj:=\langle u^{\cj-1}, u^\cj \rangle$ with
$\cj=1,\dots,\kcr+1$. (Except $u^0=(1,0)$ and $u^\kcr=(-q,n)$, their generators
$u^\cj$
are primitive lattice points (i.e.\ $\in N$) contained in
$\mbox{int}\Delta\subseteq \sigma$.)
\vspace{1ex}\\
{\em Step 1: \quad
For each $\tau^\cj$ there are $w\in E, \dalpha\in\N$ such that
$\langle u^{\cj-1}, w \rangle = \langle u^\cj, w \rangle = \dalpha$.}\\
First, it is very clear that there are a primitive lattice point $w\in M$ and
a non-negative number $\dalpha\in\R_{\geq 0}$ admitting the desired properties.
Moreover, since $u^\cj\in N$, $\dalpha$ has to be an integer, and
Reid's Lemma \zitat{P}{1}
tells us that $w\in \sigma^{\veee}$. It remains to show that $w$ belongs even
to
the Hilbert basis $E\subseteq \sigma^{\veee}\cap M$.\\
Denote by $\ell$ the length of the line segment $\overline{u^{\cj-1}u^\cj}$ on
the
``roof'' line $[\langle \bullet,\,w\rangle =\dalpha]$.
% of the cone $\tau^\cj$.
Since $\tau^\cj$ represents a T-singularity, we know from (7.3) of \cite{Homog}
(cf.\ \zitat{I}{1} of the present paper) that $\dalpha|\ell$. In particular,
$\overline{u^{\cj-1}u^\cj}$ contains the $\dalpha$-th multiple $\dalpha\cdot u$
of some lattice point $u\in\tau^\cj\cap M$ (w.l.o.g.\ not belonging to the
boundary
of $\sigma$).
Hence, $\langle u,\, w\rangle =1$ and $u\in \mbox{int}\,\sigma\cap M$,
and this implies $w\in E$.\
\vspace{1ex}\\
{\em Step 2:}
Knowing that each of the cones $\tau^1,\dots,\tau^{\kcr+1}\in\Sigma$ is
assigned
to some element $w\in E$, a
slight adaption of the notation (a renumbering) seems to be very useful:
Let $\tau^\ai=\langle u^{\ai-1},\, u^\ai\rangle$ be the cone assigned to
$w^\ai\in E$, and denote by $\dalpha_\ai,\, \ell_\ai$ the height and the length
of its ``roof'' $\overline{u^{\ai-1}u^\ai}$, respectively. Some of these cones
% (e.g.\ $\tau^1$ and $\tau^e$)
might be degenerated, i.e.\ $\ell_\ai=0$. This it at least true for the
extremal
$\tau^1$ and $\tau^e$ coinciding with the two rays spanning $\sigma$. Here we
have
even $\dalpha_1=\dalpha_e=0$; in particular $u^0=u^1=(1,0)$ and
$u^{e-1}=u^e=(-q,n)$.\\
Since $\dalpha_\ai|\ell_\ai$, we may introduce integers $k_\ai\leq a_\ai$
yielding
$\ell_\ai=(a_\ai-k_\ai)\,\dalpha_\ai$. For $\ai=2,\dots,e-1$  they are even
uniquely
determined.
\vspace{1ex}\\
{\em Step 3:}
Using the following three ingrediences
\begin{itemize}
\item[(i)]
$\langle u^{\ai-1}, \, w^\ai\rangle =
\langle u^\ai, \, w^\ai\rangle = \dalpha_\ai\quad
(\ai=1,\dots,e)\,$,
\item[(ii)]
$w^{\ai-1} + w^{\ai+1} = a_\ai\, w^\ai\quad
(\ai=2,\dots,e-1;\;$ cf.\ \zitat{CQS}{3}), and
\item[(iii)]
$\langle u^\ai - u^{\ai-1},\, w^{\ai-1} \rangle = \ell_\ai =
(a_\ai-k_\ai)\,\dalpha_\ai$ (since $\{w^{\ai-1},w^\ai\}$ forms a $\Z$-basis of
$M$),
\end{itemize}
we obtain
\[
\begin{array}{rcl}
\dalpha_{\ai-1} + \dalpha_{\ai+1}
&=&
(a_\ai\,\dalpha_\ai + \dalpha_{\ai-1}) -
(a_\ai\,\dalpha_\ai - \dalpha_{\ai+1})\\
&=&
(a_\ai\,\dalpha_\ai + \langle u^{\ai-1}, w^{\ai-1} \rangle) -
\langle u^\ai, \, a_\ai\, w^\ai- w^{\ai+1} \rangle\\
&=&
a_\ai\,\dalpha_\ai + \langle u^{\ai-1}, w^{\ai-1} \rangle
- \langle u^\ai, w^{\ai-1} \rangle\\
&=&
a_\ai\,\dalpha_\ai + \langle u^{\ai-1}-u^\ai, \,w^{\ai-1} \rangle\\
&=&
a_\ai\,\dalpha_\ai - (a_\ai-k_\ai)\,\dalpha_\ai
\; = \;
k_\ai\,\dalpha_\ai\quad (\mbox{for } \ai=2,\dots,e-1)\,.
\end{array}
\]
In particular, $k_\ai\geq 0$ (and even $\geq 1$ for $e>3$).
Moreover, since $\{w^{\ai-1},w^\ai\}$ forms
a basis of $M$ and $u^{\ai-1}\in N$ is primitive,
we have $\mbox{gcd}(\dalpha_{\ai-1},\dalpha_\ai)=1$. Hence,
$\dalpha_\ai=\qalpha_\ai$ (both series of integers satisfy the second
of the three properties mentioned in the beginning of \zitat{P}{2}).
Finally, the third of these properties yields
$[k_2,\dots,k_{e-1}]=\qalpha_1/\qalpha_2 = \dalpha_1/\dalpha_2 = 0$,
i.e.\ $\ku{k}\in K_{e-2}$.
\vspace{1ex}\\
The reversed direction (i.e.\ the fact that each $K(Y_\sigma)$-element
indeed yields a P-resolution) follows from the above calculations in
a similar manner.
\hfill$\Box$\\
\par

{\bf Remark:}
Subdividing each $\tau^\ai$ further into $(a_\ai-k_\ai)$ equal cones
(with ``roof'' length $\qalpha_\ai$ each) yields the so-called
M-resolution (cf.\ \cite{M}) assigned to a P-resolution. It is defined
to contain only T$_0$-singularities (i.e.\ T-singularities with
Milnor number 0); in exchange, $K_{\tilde{Y}|Y}$ does not need to be
relatively ample anymore. This property is replaced by ``relatively nef''.\\
\par

{\bf Examples:}
(1) The continued fraction $[1,2,2,\dots,2,1]=0$ yields $\qalpha_1=\qalpha_e=0$
and $\qalpha_\ai=1$ otherwise. In particular, the ``roof'' lines equal
$[\langle \bullet,\,w^\ai\rangle =1]$ (for $\ai=2,\dots,e-1$)
describing the RDP-resolution of $Y_\sigma$. The assigned M-resolution
equals the minimal resolution mentioned at the end of \zitat{CQS}{3}.
\vspace{1ex}\\
(2) Let us return to Example \zitat{MR}{3}:
The embedding dimension $e$ of $Y_\sigma$ is $6$, the vector
$(a_2,\dots,a_{e-1})$ equals $(2,3,2,3)$, and, except the trivial RDP element
mentioned in (1),
$K(Y_\sigma)$ contains only $(1,3,1,2)$ and $(2,2,1,3)$.\\
In both cases we already know that $\qalpha_1= \qalpha_6=0$ and
$\qalpha_2=\qalpha_5=1$. The remaining values are given by the equation
$\qalpha_3/\qalpha_4=[k_4,k_5]$, i.e.\ we obtain
$\qalpha_3=1$, $\qalpha_4=2$ or $\qalpha_3=2$, $\qalpha_4=3$, respectively.
\vspace{0.5ex}\\
Hence, in case of $(1,3,1,2)$ the fan $\Sigma$ is given by the additional
rays through $(0,1)$ and $(-4,11)$. For $\ku{k}=(2,2,1,3)$ we need the only
one through $(-1,4)$.\\
\par

%%%%%%%%%%%%%%%
%
% BIBLIOGRAPHY
%
%%%%%%%%%%%%%%%%%%

%\newpage

\end{document}